\documentstyle[amstex,12 pt]{article} \hfuzz=1.4pt %
\newtheorem{Th}{Theorem}[section]
\newtheorem{Prop}[Th]{Proposition}
\newtheorem{Lem}[Th]{Lemma}
\newtheorem{Cor}[Th]{Corollary}
\newtheorem{Rem}[Th]{Remark}
\newtheorem{Com}[Th]{Comparison with the Mabuchi-Semmes-Donaldson metric}

\newtheorem{Def}[Th]{Definition}
\newenvironment{pf}{\noindent{\it Proof. }}{$\square$\par\medskip}

\newcommand{\ddt}{\frac{d}{dt}}

\newcommand{\dmu}{\frac{\omega_\phi^n}{n!}}

\newcommand{\Rm}{{\mathbb R}}
\newcommand{\Lb}{{\mathcal L}}

\newcommand{\Hb}{{\mathcal H}}
\newcommand{\Vb}{{\mathcal V}}

\newcommand{\Cb}{{\mathcal C}}

\newcommand{\barr}{\partial}
\newcommand{\debarr}{\overline{\partial}}

\newcommand{\unmezzo}{\frac{1}{2}}

\newcommand{\arccot}{\operatorname{arccot}}

\newcommand{\Diam}{\operatorname{Diam}}

\newcommand{\Ca}{\operatorname{Ca}}
\newcommand{\Vol}{\operatorname{Vol}}
\newcommand{\Cal}{\operatorname{Cal}}
\newcommand{\sgn}{\operatorname{sgn}}
\newcommand{\dmuz}{\frac{\omega^n}{n!}}

\def\Ho{\vbox{\offinterlineskip\hbox{\kern 3pt$\scriptstyle\circ$}
\kern 1pt\hbox{$H$}}}
\begin{document}
\title{The Calabi's metric for the space of K\"{a}hler metrics}
\author{Simone CALAMAI}

\date{}

\maketitle
\begin{abstract}
\noindent Given any closed K\"{a}hler manifold we define,
following an idea by Eugenio Calabi \cite{CC2}, a Riemannian
metric on the space of K\"{a}hler metrics regarded as a infinite
dimensional manifold. We prove several geometrical features of the
resulting space, some of which we think were already known to
Calabi. In particular, the space has positive constant sectional
curvature and admits explicit unique smooth solutions for the
Cauchy and the Dirichlet problems for the geodesic equation.
\end{abstract}

\section{Introduction and main results}

\noindent The space of K\"{a}hler metrics is the natural
environment to study the problem about the existence of Extremal
K\"{a}hler metrics, which goes back to E. Calabi \cite{Cal} and is
a central question in K\"{a}hler geometry. The space of K\"{a}hler
metrics is also the setting for the K\"{a}hler-Ricci flow and the
Calabi flow. Suppose $(M, \omega)$ is a closed K\"{a}hler manifold
with K\"{a}hler form $\omega$. In this introduction and in the
remainder we refer to the space of K\"{a}hler metrics with the
following notation
$$
\Cb:=\{ u \in C^{\infty}(M , \Rm) \, | \, \int_M e^u
\frac{\omega^n}{n!} = \int_M \frac{\omega^n}{n!}\}.
$$
\noindent This work is the proposal of a geometric structure for
the space of K\"{a}hler metrics viewed as an infinite dimensional
manifold. Another geometry for the same space turned out to be
decisive to affirmatively answer the question of uniqueness of
Extremal metrics. For that, great contributions were given by Xiu
Xiong Chen \cite{Chen} and S. Donaldson \cite{Don}. In that case,
the geometry arises from a metric which we call the
Mabuchi-Semmes-Donaldson's metric, since those three authors
defined it and contributed to its study. Our geometry arises from
an intuition by E. Calabi, as we learned by Xiu Xiong Chen
\cite{CC2}. The expression of the Calabi's metric is
$$
<v , w>_u = \int_M vwe^u \frac{\omega^n}{n!}, \qquad v,w \in T_u
\Cb.
$$
The
purpose of this paper is to compare those two geometries;
hereafter we summarize the main results. For the first one,
compare Theorems \ref{teolevicivitaexistence},
\ref{teoseccurvcostpos} of the present paper.
\begin{Th}
Let $(M, \omega)$ a closed K\"{a}hler manifold. The Calabi's
metric admits the Levi Civita covariant derivative; its sectional
curvature is positive, constant and equal to $s= \frac{1}{4\Vol}$,
where $\Vol$ is the volume of the manifold $M$. As a consequence,
the space of K\"{a}hler metrics is a locally symmetric space.
\end{Th}
\noindent In his seminal paper, T. Mabuchi \cite{Mab} got similar
results for his metric; i.e. he explicitly computed the Levi
Civita covariant derivative which he proved to entail a structure
of locally symmetric space. But he proved also that his metric has
non positive sectional curvature. On this aspect, that metric
seems to be opposite to the Calabi's one, and of course the two
models are not isometric. The Calabi's metric case already showed
the pleasant feature to have constant curvature; more pleasant
features comes with the analysis of the geodesic equation (cfr.
Theorem \ref{teocauchyproblemgeodesics}).
\begin{Th}
The geodesic equation is equivalent to an ordinary differential
equation, namely to the equation
$$
(e^{\frac{u}{2}})_{tt} + e^{\frac{u}{2}}=0.
$$
\end{Th}
\noindent On this side, our case is the best one can hope for and
again very different from the other one. There, unaware f the
Mabuchi's paper, S. Semmes \cite{Sem} rediscovered the same metric
via a totally different line of reasoning. Indeed, he was studying
the Homogeneous Complex Monge-Amp\`ere equation on a complex
domain $D\subset {\mathbb{C}}^n$ and proved that in some special
cases it is a geodesic equation for the space of plurisubharmonic
functions on $D$. Moreover, he was able to construct a metric on
the space of plurisubharmonic functions on $D$ such that its
corresponding geodesic equation is precisely the Homogeneous
Complex Monge-Amp\`ere equation. In \cite{Don}, S. Donaldson made
precise that the geodesic equation for the
Mabuchi-Semmes-Donaldson's metric is equivalent to a Homogeneous
Complex Monge-Amp\`ere equation. What really interests us is not
the complexity of that equation, but rather the consequences on
the geometry; these appear in the next result (cfr. Theorems
\ref{teocauchyproblemgeodesics},\ref{teoexpinjandsurj} and
\ref{teonoconjpoints}).
\begin{Th}
The Cauchy problem for the geodesic equation has unique real
analytic solutions for any initial data; the explicit expression
of a solution is
$$
u(t) = u_0 + 2\log (\cos (t) + \frac{v_0}{2}\sin (t)).
$$
Moreover, this allows to define the exponential map, which is
injective and surjective. Finally, the space of K\"{a}hler metrics
does not have conjugate points.
\end{Th}
\noindent In \cite{Don}, S. Donaldson showed examples of
nonexistence of solutions for the Cauchy problem in the
Mabuchi-Semmes-Donaldson's setting. Thus, the notions of
exponential map and conjugate points couldn't be carried over. On
this aspect, we may say that the Calabi's geometry is richer than
that one. This appears also in the next result (cfr. Theorems
\ref{teogeodesicallyconvex}, \ref{disadistance},
\ref{teodiameter}).
\begin{Th}
The Dirichlet problem for the geodesic equation has unique real
analytic solutions for any boundary data $u_0, u_1 \in \Cb$; the
explicit expression of a solution is
$$
e^{\frac{u(t)}{2}} = e^{\frac{u_0}{2}}(\cos(t)+ e^{\frac{u_1 -
u_0}{2}}\frac{\sin(t)}{\sin(t_0)} -
\cos(t_0)\frac{\sin(t)}{\sin(t_0)}).
$$
The space of K\"{a}hler metrics is a metric space with distance
function real analytic. Geodesics are minima of the length
functional. Moreover, the diameter of the space is
$\frac{\pi}{2}R$, where the sectional curvature is considered as
$s=\frac{1}{R^2}$.
\end{Th}
\noindent In \cite{Chen}, Xiu Xiong Chen solved the Dirichlet
problem for geodesics in the Mabuchi-Semmes-Donaldson's setting.
Those geodesic segments, called Chen's geodesics, are of class
$C^{1,1}$ and they stay in the closure of the space of K\"{a}hler
metrics. The improvement of the latter two facts is still an open
problem. Moreover, Chen showed that the space is a metric space
with distance function of class $C^1$, and that geodesics minimize
the length functional. We mention that these results led to prove
the uniqueness of constant scalar curvature metrics; this was done
by Xiu Xiong Chen \cite{Chen} and by Xiu Xiong Chen and Gang Tian
\cite{CT}. Other important results are: E. Calabi and Xiu Xiong
Chen \cite{CC} proved that the space of K\"{a}hler potentials
endowed with the Mabuchi-Semmes-Donaldson's metric is non
positively curved in the sense of Alexandrov and Xiu Xiong Chen
\cite{ChenIII} proved the inequality involving the geodesic
distance induced by the Mabuchi-Semmes-Donaldson's metric $
E(\phi_1) - E(\phi_0) \leq \sqrt{Cal (\phi_1)} d(\phi_0 , \phi_1).
$ Recently, Xiu Xiong Chen and Song Sun \cite{ChenSun} reproved
the latter two results using a quantization of the space of
K\"{a}hler potentials. About the question on the diameter, the two
situations seems to be opposite; in fact Donaldson showed how the
study of geodesic rays in the space of K\"{a}hler potentials
endowed with the Mabuchi-Semmes-Donaldson's metric that would lead
to results about the existence of Extremal metrics (on this
aspect, see \cite{AT} and \cite{ChenIV}). In the Calabi's case,
the fact that the sectional curvature is constant  and positive,
the fact that there are no conjugate points and the explicit value
of the diameter led to the result which summarize the present
study of the geometry of the Calabi's metric. (cfr. Theorems
\ref{teoimmersionsphere}, \ref{teodistancefromtheboundaryiszero}).
\begin{Th}
There is an isometric immersion of the space of K\"{a}hler metrics
into the space of real valued smooth functions on $M$ endowed with
the Euclidean metric, and the image of the space $\Cb$ under this
immersion is a portion of an infinite dimensional sphere.
Moreover, the distance to the boundary of any point of the space
is zero.
\end{Th}

\smallskip

\noindent Recently the interest for Riemannian metrics in infinite
dimensional spaces has been renewed by S. Donaldson \cite{Don2}
who generalized the considerations made in the K\"{a}hler case to
the Riemannian case, i.e. when is given a closed Riemannian
manifold $(M,g)$ neither necessarily K\"{a}hler nor complex. In
this more general case the Riemann metric is defined no more in a
fixed K\"{a}hler class but in the space of volume forms which are
conformal to the given one $d\mu$ and which have the same area;
$$
V := \{ f \in C^{\infty}(M, (0,+\infty)) \; | \; \int_M f d\mu =
\int_M d\mu \}.
$$
The  space $V$ can be parameterized as
$$
{\mathcal S} := \{ \phi \in C^{\infty}(M, \Rm) \; | \; 1+ \Delta_g
\phi
>0 \}
$$
and the metric in ${\mathcal S}$ is defined by
$$
<\psi , \chi >_\phi := \int_M \psi \chi (1 + \Delta_g \phi ) d\mu
.
$$
When $(M, g)$ has real dimension $2$ this metric coincides with
the Mabuchi-Semmes-Donaldson's metric. Endowed with this metric by
Donaldson, $V$ has a very interesting geometry structure as well.
The existence and the study of the geometry arising from
Donaldson's metric on $V$ is due to Xiu Xiong Chen and Weiyong He
(see \cite{CH} for more details), and their arguments generalize
the techniques used by Chen in \cite{Chen}.\\
\noindent We remark here that we can generalize the Calabi's
metric to the space $V$; in fact we parameterize the space $V$ as
$$
\Cb := \{u \in C^{\infty}(M, \Rm) \;|\; \int_M e^u d\mu = \int_M
d\mu \}
$$
and we define the metric
$$
<v,w>_u := \int_M vw e^u d\mu.
$$
We claim that the above definition coincides with the Calabi's
metric. In fact, in the K\"{a}hler case, the Calabi's volume
conjecture \cite{Cal0} solved by S.T. Yau \cite{Yau} entails  that
the space of metrics in the anti-canonical bundle of the given
K\"{a}hler manifold is equivalent to the space of K\"{a}hler
potentials up to a constant; moreover it is not hard to see that
through the above equivalence the Calabi's metric translates into
a $L^2$ metric on the space of metrics in the anti-canonical
bundle of the given K\"{a}hler manifold. \\
\noindent An interesting reference for a further generalization of
the Calabi's metric are the recent lecture notes by Bourguignon
\cite{B}, which relates together spaces of probability theory,
optimal transportation theory and Riemannian geometry.\\
\noindent \textbf{Acknowledgements:} I want to thank my advisor
Giorgio Patrizio, Claudio Arezzo and Xiu Xiong Chen for his
support. Thanks also to Chen's group of students, in particular to
Song Sun for many helpful discussions and his interest on this
problem. Part of this work was written at UW-Madison and at USTC,
Hefei. I thank both those Universities for their hospitality.

\section{The Calabi's metric}
In this section we introduce the Calabi's metric. For the reader's
convenience, we recall some well known definitions which lead to
that notion. From now on, $(M, \omega)$ is a closed K\"{a}hler
manifold of complex dimension $n$ with K\"{a}hler form $\omega$.

\begin{Def}
\emph{The space of volume conformal factors of $(M , \omega)$ is
$$
\Cb := \{ u \in C^{\infty}(M , \Rm) \, | \, \int_M e^{u} \dmuz =
\int_M \dmuz \}.
$$
The space of K\"{a}hler potentials $\Hb$ is $ \Hb := \{ \phi \in
C^{\infty} (M , {\mathbb R}) \, | \, \omega + i \barr \debarr \phi
>0\}. $ The space $\tilde{\Hb}$ of normalized potentials is given
by $ \tilde{\Hb} := \{ \phi \in \Hb \, | \, L(0,\phi ) = 0\}, $
where $L(\rho , \eta):= \frac{1}{Vol} \int_0^1 (\int_M \phi
\frac{\omega_{s\phi}^n}{n!}) ds \mid_{\phi = \rho}^{\phi = \eta}$
and $\rho,\eta \in \Hb$, according to \cite{Mab} and \cite{Mab2}.
}
\end{Def}

\begin{Rem}\label{calmap}
\emph{The space $\Cb$ is diffeomorphic to the space $\tilde{\Hb}$;
indeed the map $\Cal : \tilde{\Hb} \rightarrow \Cb$, $\Cal (\phi)
= \log \left( \frac{\omega_\phi^n}{\omega^n}\right) $ is a
diffeomorphism (see \cite{Bes}) according to the Calabi's
conjecture \cite{Cal0} solved by Yau \cite{Yau}.}
\end{Rem}

\begin{Def}
\emph{ A smooth curve in $\Cb$ is a map $u= u(t): (-\epsilon ,
\epsilon) \rightarrow \Cb$ such that the map
$$
\begin{array}{cccc}
  u: & (-\epsilon , \epsilon)\times M & \rightarrow  &\Rm\\
   & (t,p) & \mapsto & (u(t))(p). \\
\end{array}
$$}
is smooth.
\end{Def}

\begin{Def}
\emph{ Fix a point $u \in \Cb$. Let $\alpha : (-\epsilon ,
\epsilon) \rightarrow \Cb$ a smooth curve with $\alpha (0)= u$.
Then $\frac{d\alpha}{dt}_{|t=0}$ is a tangent vector at $u$. The
whole tangent vectors at $u$ form $T_u \Cb$, the tangent space at
$u$. A similar definition holds for $\Hb$ and $\tilde{\Hb}$.}
\end{Def}

\begin{Prop}
The following characterization of the tangent space at $u$ holds
$$
T_u \Cb = \{ v \in C^{\infty}(M  ,\Rm) \, | \, \int_M ve^u \dmuz =
0\}.
$$
Similarly, $T_\phi \tilde{\Hb} = \{ \psi \in C^{\infty}(M , \Rm)
\, | \, \int_M \psi \dmu =0\}$.
\end{Prop}

\begin{Def}
\emph{ Let $u= u(t):(-\epsilon , \epsilon)\rightarrow \Cb$ be a
smooth curve; a real function $v\in C^{\infty}((-\epsilon ,
\epsilon)\times M , \Rm)$ is a smooth section on $u$ when
$v(t,\cdot)\in T_{u(t)}\Cb$ for any $t$.}
\end{Def}

\begin{Def}\label{generaldefofriemannianmetric}
\emph{ A metric or 'Riemannian metric' on $\Cb$ is a positive
bilinear form at any tangent space of $\Cb$, which is
differentiable along any smooth sections on any smooth curve on
$\Cb$.}
\end{Def}

\begin{Def}
\emph{ The Calabi's metric is given by, at any $u\in \Cb$,
$$
\Ca(v,w) = <v , w>_u := \int_M vwe^u \dmuz , \quad v,w\in T_u \Cb.
$$}
\end{Def}

\begin{Rem}
\emph{ The Calabi's metric is an idea by Eugenio Calabi \cite{CC2}
and we learned it from Xiu Xiong Chen. It satisfies the
requirements of the Definition
\ref{generaldefofriemannianmetric}.}
\end{Rem}

\begin{Com}
\emph{ The map $\Cal$ recalled in Remark \ref{calmap} is an
isometry between $(\Cb , \Ca)$ and the space $\tilde{\Hb}$ endowed
with the metric
$$
<\psi , \chi >_{\phi} = \int_M (\Delta_\phi \psi)(\Delta_\phi
\chi) \dmu, \quad \psi, \chi \in T_{\phi} \tilde{\Hb};
$$
the Mabuchi-Semmes-Donaldson metric is (isometric to)
$$
\prec \psi , \chi \succ_{\phi} := \int_M \psi \chi \dmu, \quad
\psi, \chi \in T_{\phi} \tilde{\Hb}.
$$
So the Calabi's metric is a $L^2$ pairing of the Laplacian instead
of just a $L^2$ pairing. }
\end{Com}

\section{The Levi Civita covariant derivative}

\noindent In this section we compute the  explicit expression of
the Levi Civita covariant derivative for the Calabi's metric. This
a fortiori proves its existence which in this infinite dimensional
setting is not guaranteed by the standard finite dimensional
argument (cfr. \cite{Don}).

\begin{Def}\label{levicivitacovariantderivative}
\emph{The Levi Civita covariant derivative for the space $\Cb$
endowed with a 'Riemannian metric' $g$ is a map defined on every
smooth curve $u = u(t) : (-\epsilon , \epsilon ) \rightarrow \Cb$
and a smooth real section $v$ along $u$. The Levi Civita covariant
derivative of $v$ along $u$ is denoted by $D_t v$ and it is a
smooth section along $u$. It is required to satisfy
$$
\begin{array}{l}
  (i) D_t (v+w ) = D_t v + D_t w ;\\
  (ii) D_t (fv ) = f D_t v + \ddt f v;\\
  (iii) \ddt g_\phi( v , w ) = g_\phi (D_t v , w ) + g_u (v , D_t w) ;\\
  (iv) \tau (\alpha ) : =  D_t \frac{\barr \alpha}{\barr s} - D_s \frac{\barr \alpha}{\barr t} =0,\\
\end{array}
$$
where, along the path $u$, are given $v,w $ smooth sections and
$f$ is a smooth function; $\alpha$  is a smooth two parameter
family in $\Cb$ and $\tau$ is the torsion.}
\end{Def}

\begin{Rem}
\emph{The requirement for the Levi Civita covariant derivative to
be torsion free finds an application in those computations where
there is a smooth two parameter family in the space of K\"{a}hler
metrics; we are allowed to switch in this sense $
 D_t  \frac{\barr \alpha }{\barr s } = D_s \frac{\barr \alpha }{\barr t
 }
$, as follows from the fact that $D$ is torsion free.}
\end{Rem}

\begin{Rem}
\emph{We used  the expression '\emph{the} Levi Civita covariant
derivative'. In fact it holds both in the finite dimensional and
in this infinite dimensional case that the Levi Civita covariant
derivative, if it exists, is unique (cfr. \cite{Don}). For a
finite dimensional Riemannian manifold the Fundamental lemma of
Riemannian geometry (see \cite{Mil}) guarantees the existence of
the Levi Civita covariant derivative; that argument does not work
in our infinite dimensional setting (cfr. \cite{Don}). This
motivates the next result.}
\end{Rem}

\begin{Th}\label{teolevicivitaexistence}
The Levi Civita covariant derivative for the space $\Cb$ endowed
with the Calabi's metric exists; moreover its explicit expression
is
$$
D_t v = v' + \unmezzo vu' + \frac{1}{2Vol}\int_M vu' e^u \dmuz ,
$$
where $u = u(t) :(-\epsilon , \epsilon) \rightarrow \Cb$ is a
smooth path and $v=v(t)$ is a smooth section along $u$.
\end{Th}

\smallskip

\noindent \begin{pf} It is sufficient to prove the second part of
the statement. Let $u, v,w, f$ be as in Definition
\ref{levicivitacovariantderivative}. To prove that $D_t v$ is a
smooth section along $u$ it is enough to show that $\int_M (D_t v)
e^u \dmuz =0$;
$$
\int_M (\ddt v +\unmezzo \ddt u v +\frac{1}{2Vol}\int_M v \ddt u
\dmuz) e^u \dmuz =
$$
$$
= \int_M \ddt v e^u \dmuz + \int_M v \ddt u e^u \dmuz = \ddt
\int_M v e^u \dmuz  = 0,
$$
as claimed. The following computation
$$
D_t (v + w) = \ddt (v + w) + \unmezzo (v+w)\ddt u +
\frac{1}{2Vol}\int_M (v+w)\ddt u \dmuz =
$$
$$
=(\ddt v  + \unmezzo v\ddt u + \frac{1}{2Vol}\int_M v \ddt u
\dmuz)+(\ddt w  + \unmezzo w\ddt u + \frac{1}{2Vol}\int_M w \ddt u
\dmuz) = D_t v + D_t w
$$
proves $(i)$ of Definition \ref{levicivitacovariantderivative}.
About $(ii)$,
$$
D_t (fv) = \ddt (fv) + \unmezzo (fv) \ddt u +\frac{1}{2Vol}\int_M
(fv) \ddt u \dmuz =
$$
$$
= v\ddt f + (f\ddt v + f\unmezzo v\ddt u +f \frac{1}{2Vol}\int_M v
\ddt u \dmuz ) = v\ddt f + f D_t v ,
$$
as required, where is used that $f$ does not depend on spacial
variables. To check $(iii)$ with $g$ being the Calabi's metric it
is enough to show
$$
\ddt \int_M v^2 e^u \dmuz = 2\int_M v(D_t v) e^u \dmuz .
$$
The above equality holds since
$$
\ddt \int_M v^2 e^u \dmuz = \int_M (2v\ddt v + v^2 \ddt u)e^u
\dmuz =
$$
$$
= 2\int_M  v (\ddt v + \unmezzo v \ddt u  + \frac{1}{2Vol}\int_M )
e^u \dmuz - \int_M v \frac{1}{Vol} <v , \ddt u>_u e^u \dmuz =
$$
$$
= 2 \int_M v (D_t v) e^u \dmuz - \frac{1}{Vol} <v , \ddt u >_u
\int_M v e^u \dmuz = 2 \int_M v(D_t v ) e^u \dmuz - 0.
$$
Finally, about $(iv)$, let $\alpha = \alpha (s,t)$ be a smooth two
parameter family. The computation
$$
D_t \frac{\barr \alpha}{\barr s} = \frac{\barr^2 \alpha}{\barr t
\barr s} + \unmezzo \frac{\barr \alpha}{\barr s}\frac{\barr
\alpha}{\barr t}  + \frac{1}{2Vol}<\frac{\barr \alpha}{\barr s},
\frac{\barr \alpha}{\barr t}>_\alpha = D_s \frac{\barr
\alpha}{\barr t}
$$
completes the proof. \end{pf}

\section{The Calabi's metric has positive constant\\ sectional curvature}

In \cite{CC2} Eugenio Calabi was aware of the fact that his metric
has positive curvature. In this section we prove that the space
$\Cb$ endowed with the Calabi's metric has sectional curvature
equal to $(4\Vol)^{-1}$, where $\Vol=\int_M \dmuz$ is the volume
of $M$. As a consequence, the space is locally symmetric, that is
the covariant derivative of the curvature tensor is identically
zero, and the Calabi's metric is not isometric to the
Mabuchi-Semmes-Donaldson metric. For the sake of notations, we
will write $u_s $ for $\frac{\barr u}{\barr s}$ and also for
$\frac{\barr u}{\barr s}_{|s=0}$ when no confusion arises;
moreover we will write $<\cdot , \cdot>$ for $<\cdot , \cdot>_u$,
i.e. with the omission of $u$.

\begin{Def}
\emph{Let $g$ be a metric on $\Cb$ with Levi Civita covariant
derivative $D$; the curvature tensor is defined on a four
parameter family $u (q ,r, s, t)$ by
$$
R(u_q , u_r , u_s , u_t ) := g(R(u_q , u_r ) u_s , u_t ),
$$
where $R(u_q , u_r ) u_s := (D_q D_r - D_r D_q) u_s$. Its
covariant derivative is a $5$-tensor given by
$$
D R(u_p , u_q , u_r, u_s ,u_t ) := D_p R(u_q , u_r, u_s , u_t) -
R(D_p u_q , u_r, u_s ,u_t)- \cdots - R( u_q , u_r, u_s , D_p u_t).
$$
The sectional curvature of the plane $ P:= \Rm \frac{\barr
u}{\barr s} + \Rm \frac{\barr u}{\barr t}$ is given by
$$
K(P) := -\frac{g(R(u_s,u_t)u_s ,u_t)}{g(u_t,u_t) g(u_s , u_s) -
g(u_s ,u_t)^2}.
$$}
\end{Def}

\begin{Th}\label{curvaturetensorcputation}
Let $u = u(q,r,s,t)$ a four parameter family in $\Cb$; $u :
(-\epsilon, \epsilon )^4 \rightarrow \Cb$. The curvature tensor
for the Calabi's metric is
$$
R(u_q , u_r , u_s , u_t ) = \frac{1}{4\Vol}( < u_r , u_s>_u < u_q
,u_t
>_u - <u_q  ,u_s
>_u < u_r , u_t >_u).
$$
\end{Th}

\smallskip

\noindent {\it Proof.} It is convenient to compute first $D_q D_r
u_s $. By definition
$$
D_q D_r u_s = D_q (u_{sr} + \unmezzo u_s u_r + \frac{1}{2\Vol}<u_s
, u_r >)=
$$
$$
= u_{srq} + \unmezzo u_{sr} u_q + \frac{1}{2\Vol}<u_{sr} , u_q>
+\unmezzo (u_s u_r )_q + \frac{1}{4}u_s u_r u_q +
\frac{1}{4\Vol}<u_s u_r , u_q> +
$$
$$
+ \frac{1}{2\Vol}\frac{\barr}{\barr q}< u_s , u_r > +
\frac{1}{4\Vol}<u_s , u_r> u_q + \frac{1}{4\Vol^2} <<u_s , u_r> ,
u_q >=
$$
$$
=u_{srq} + \unmezzo (u_{sr} u_q + u_{sq} u_{r}) + \unmezzo u_s
u_{rq} +\frac{1}{4} u_s u_r u_q +  \frac{1}{2\Vol} <u_{sr},u_q>+
$$
$$
+\frac{1}{4\Vol} <u_s u_r ,u_q> + \frac{1}{2\Vol}<u_{sq} +
\unmezzo u_s u_q + \frac{1}{2\Vol} <u_s , u_q >  , u_r > +
$$
$$
+\frac{1}{2\Vol}<u_s , u_{qr} + \unmezzo  u_r u_q +
\frac{1}{2\Vol} <u_r , u_q >  > + \frac{1}{4\Vol}<u_s , u_r > u_q
=
$$
$$
=u_{srq} + \unmezzo (u_{sr} u_q + u_{sq} u_{r}) + \unmezzo u_s
u_{rq} +\frac{1}{4} u_s u_r u_q +
$$
$$
+  \frac{1}{2\Vol}( <u_{sr},u_q> + <u_{sq} , u_r > ) +
\frac{1}{4\Vol} (<u_s u_r , u_q > + <u_s u_q , u_r >) +
$$
$$
+\frac{1}{2\Vol}<u_s , u_{qr}> + \frac{1}{2\Vol}<u_s , u_q u_r > +
\frac{1}{4\Vol}<u_r , u_s > u_q .
$$
The above computation leads to this explicit expression
$$
R(u_q , u_r) u_s = D_q D_r u_s - D_r D_q u_s =
\frac{1}{4\Vol}(<u_r , u_s
> u_q - <u_q , u_s > u_r),
$$
form which the claim of the theorem easily follows.$\Box$

\noindent This first consequence of the Theorem
\ref{curvaturetensorcputation} is immediate.
\begin{Cor}
The curvature tensor for the Calabi's metric is locally symmetric,
i.e. $D R =0.\;\Box$
\end{Cor}

\noindent The following consequence of the previous theorem is so
important that we label  again as a theorem its statement.

\begin{Th} \label{teoseccurvcostpos}
The curvature of the space $\Cb$ endowed with the Calabi's metric
is $\frac{1}{4\Vol}$, i.e. it is positive, constant and it depends
only on the volume of the manifold $(M , \omega)\;\Box$.
\end{Th}

\begin{Com} \emph{ Both metrics induce a structure of locally symmetric
space, i.e. their curvature tensor is covariant constant (see
\cite{Don}). Mabuchi in \cite{Mab} proved his metric to have non
positive sectional curvature;  Theorem \ref{teoseccurvcostpos}
entails that the two metrics are not isometric.}
\end{Com}

\begin{Rem} \emph{We will often employ the notation $R^2 = \frac{1}{4\Vol}$; the normalization $R=1$ corresponds
to replace $\omega$ with $\frac{1}{\sqrt[n]{4\Vol}}\omega$.}
\end{Rem}

\section{The Cauchy problem for the geodesic\\ equation}

In this section we prove that the Cauchy problem for the geodesic
equation admits, for any initial datum, an explicit unique smooth
solution; as remarked by S. Donaldson in \cite{Don}, this does not
happen for the Mabuchi-Semmes-Donaldson metric. We define the
exponential map and we prove that it is a diffeomorphism on its
domain. We point out that in our case the geodesic equation
reduces to a ordinary differential equation, while for the
Mabuchi-Semmes-Donaldson metric it is equivalent to a
Monge-Amp\`ere equation.

\begin{Def}
\emph{ Let $g$ be 'Riemannian metric' for the space $\Cb$. Assume
that $\nabla_t$ is the Levi Civita covariant derivative induced by
the metric $g$. A path $\phi = \phi (t)$ is said to be a geodesic
curve if it satisfies the geodesic equation
$$
\nabla_t \frac{\barr \phi }{\barr t} =0.
$$}
\end{Def}

\begin{Th}\label{teocauchyproblemgeodesics}
Let $u_0\in \Cb$ and $v_0 \in T_{u_0}\Cb$ be respectively the
given initial position and initial velocity. Then there exists a
unique geodesic curve $u = u(t)$ in $\Cb$ satisfying the Cauchy
problem
$$
\left\{
\begin{array}{c}
  D_t u'=0 \\
  u(0)=u_0 \\
  u'(0)= v_0. \\
\end{array}
\right.
$$
The curve is smooth and, if $v_0 \neq 0$, is given by
$$
u(t)= u_0 + 2\log(\cos(\frac{|v_0|_{u_0}}{R}t) + \frac{v_0}{2}
\frac{R}{|v_0|_{u_0}} \sin (\frac{|v_0|_{u_0}}{R}t));
$$
it is defined for those values of $t$ in the open interval
$$
(- \frac{R}{|v_0|_{u_0}}\arccot(\frac{\max
v_0}{2}\frac{R}{|v_0|_{u_0}}),
\frac{R}{|v_0|_{u_0}}\arccot(\frac{-\min
v_0}{2}\frac{R}{|v_0|_{u_0}})),
$$
where the function $\arccot$ is meant to range in the interval
$(0,\pi)$.
\end{Th}

\begin{pf} Fix $R=1$, so that the geodesic equation is $u''+\unmezzo (u')^2 + 2<u' , u'>=0$.
A first claim  is that along a geodesic $<u', u'>(t)=<u' ,
u'>(0)$; in fact $\ddt <u' , u'> = 2<Du' , u'> =0$ by the metric
compatibility. Introduce the notation $|v_0|_{u_0} := <u'
,u'>(0)$; the geodesic equation is thus proved to be equivalent to
the second order ordinary differential equation
$$
u'' + \frac{1}{2}(u')^2 + 2|v_0|_{u_0}^2 =0.
$$
This is in fact the harmonic oscillator equation of proper
frequency $|v_0|_{u_0}$ in the unknown $e^{\frac{u}{2}}$; indeed
$(e^{\frac{u}{2}})''=\unmezzo e^{\frac{u}{2}}(u''+\unmezzo
(u')^2)$, so that the geodesic equation is equivalent to
$$
(e^{\frac{u}{2}})'' +|v_0|_{u_0}^2 e^{\frac{u}{2}} =0 ,
$$
having multiplied it by the positive factor $\unmezzo
e^{\frac{u}{2}}$. About the corresponding Cauchy problem, in the
case that $v_0=0$ readily the unique solution is $u(t)= u_0$ for
any $t$; if $|v_0|_{u_0}>0$, the unique solution is the smooth
curve
$$
e^{\frac{u}{2}} = e^{\frac{u_0}{2}}(\cos(|v_0|_{u_0}t) +
\frac{v_0}{2|v_0|_{u_0}}\sin( |v_0|_{u_0} t)),
$$
which is equivalent to that one in the statement. Notice that the
solution makes sense as long as $(\cos(|v_0|_{u_0}t) +
\frac{v_0}{2|v_0|_{u_0}}\sin( |v_0|_{u_0} t))>0$, which gives the
interval of definition of the geodesic that is written in the
statement. The case when $R\neq 1$ follows exactly from the same
argument.
\end{pf}

\noindent Motivated by the last result we try a definition of the
exponential map.

\begin{Def}\label{defzofexpmap}
\emph{ Let $\gamma (t , u_0 , v_0) \in \Cb$ be the point reached
after a time $t$ by the geodesic with initial data $(u_0 , v_0 )$.
The exponential map at $u_0$ assigns to $v_0 \in T_{u_0} \Cb$ the
quantity $\exp_{u_0}(v_0)= \gamma(1 , u_0 , v_0)$ whenever it
makes sense, i.e. in
$$
D_{u_0} := \{ v\in T_{u_0} \Cb \, | \, |v|_{u_0} < \arccot
(\frac{-\min v}{2|v|_{u_0}}) \}.
$$
}
\end{Def}

%

\begin{Th}\label{teoexpinjandsurj}
For any $u_0 \in \Cb$, the exponential map $exp_{u_0} : D_{u_0}
\rightarrow \Cb$ is injective and surjective.
\end{Th}

\smallskip

\noindent \begin{pf} About the injectivity, first notice that
$\exp_{u_0} (v) = u_0$ if and only if $v=0$. Then suppose
$\exp_{u_0} (v) = \exp_{u_0} (w)$ with $v,w\neq 0$, i.e.
$$
\cos (|v|_{u_0}) + \frac{v}{2|v|_{u_0}}\sin(|v|_{u_0}) = \cos
(|w|_{u_0}) + \frac{w}{2|w|_{u_0}}\sin(|w|_{u_0}).
$$
An integration gives $\cos(|v|_{u_0}) = \cos(|w|_{u_0})$; since
any nonzero  $v \in D_{u_0} $ satisfies $0< |v|_{u_0} \leq
\frac{\pi}{2} $ that entails $|v|_{u_0} = |w|_{u_0}$. The
hypothesis now writes as $(v-w) \frac{\sin(|v|_{u_0})}{2|v|_{u_0}}
=0$, that gives $v=w$ as claimed. To see that $\exp_{u_0}$ is
surjective, let $w\in \Cb, w\neq u_0$ be a generic element. It is
required to prove that there is a solution $0\neq v\in D_{u_0}$
for the equation $ e^{\frac{u_0}{2}}(\cos (|v|_{u_0}) +
\frac{v}{2|v|_{u_0}}\sin (|v|_{u_0})) = e^{\frac{w}{2}}. $
Multiply by $e^{\frac{u_0}{2}}$ and integrate to get $ \cos
(|v|_{u_0}) = 4 \int_M e^{\frac{u_0 +w}{2}}\frac{\omega^n}{n!}\in
(0, 1). $ Insert the last result in the first equation to get this
expression of $v$ in terms of $u_0$ and $w$
$$
v =  (e^{\frac{w-u_0}{2}} - 4 \int_M
e^{\frac{u_0+w}{2}}\frac{\omega^n}{n!}) \frac{2 \arccos(4 \int_M
e^{\frac{u_0+w}{2}}\frac{\omega^n}{n!})}{\sin (\arccos(4 \int_M
e^{\frac{u_0+w}{2}}\frac{\omega^n}{n!}))}.
$$
It remains only to check that $v$ is in the domain of definition
of the exponential map, i.e. $ \arccot (\frac{-\min v}{2|v|_{u_0}
}) - |v|_{u_0} >0. $ Notice that
$$
-\min v =(4 \int_M  e^{\frac{u_0 +w}{2}}\frac{\omega^n}{n!}) -
\min e^{\frac{w-u_0}{2}}) \frac{2 \arccos(4 \int_M e^{\frac{u_0
+w}{2}}\frac{\omega^n}{n!})}{\sin (\arccos(4 \int_M e^{\frac{u_0
+w}{2}}\frac{\omega^n}{n!}))};
$$
thus
$$
\arccot (\frac{-\min v}{2|v|_{u_0} }) = \arccot ( \frac{4 \int_M
e^{\frac{u_0 +w}{2}}\frac{\omega^n}{n!} - \min e^{\frac{w-u_0
}{2}}}{\sin (\arccos(4 \int_M e^{\frac{u_0
+w}{2}}\frac{\omega^n}{n!}))}).
$$
Also
$$
|v|_u = \arccot (\frac{4 \int_M e^{\frac{u_0
+w}{2}}\frac{\omega^n}{n!}}{\sin(\arccos(4 \int_M e^{\frac{u_0
+w}{2}}\frac{\omega^n}{n!}))});
$$
since the function $\arccot$ is decreasing and $\frac{- \min
e^{\frac{w-u_0}{2}}}{\sin (\arccos(4 \int_M e^{\frac{u_0
+w}{2}}\frac{\omega^n}{n!}))}<0$, conclude
$$
\arccot ( \frac{4 \int_M  e^{\frac{u_0 +w}{2}}\frac{\omega^n}{n!}
- \min e^{\frac{w-u_0}{2}}}{\sin (\arccos(4 \int_M e^{\frac{u_0
+w}{2}}\frac{\omega^n}{n!}))}) - \arccot (\frac{4 \int_M
e^{\frac{u_0 +w}{2}}\frac{\omega^n}{n!}}{\sin(\arccos(4 \int_M
e^{\frac{u_0 +w}{2}}\frac{\omega^n}{n!}))}) >0,
$$
as wanted.
\end{pf}

\begin{Lem}\label{lemmahomogeneity}
Let $a\in \Rm$ and $\gamma (t, u_0 , v_0)$ be as in Definition
\ref{defzofexpmap}; then as long as both sides  make sense we have
$$
\gamma (t , u_0 , a v_0 ) = \gamma( at , u_0 , v_0 ).
$$
\end{Lem}

\smallskip

\noindent \begin{pf} If $a$ or $t$ are zero, then the formula
reads $u_0 = u_0$; so consider the case when both $t$ and $a$ are
non zero. The lefthand side reads
$$
\gamma (t , u_0 , a v_0 ) = u_0 + 2\log(\cos(|av_0 |_{u_0}t) +
\frac{av_0}{2|av_0|_{u_0}}\sin(|av_0 |_{u_0}t))=
$$
$$
=u_0 + 2\log(\cos(|a|\cdot |v_0 |_{u_0}t) +
\frac{\sgn(a)v_0}{2|v_0 |_{u_0}}\sin(|a| \cdot |v_0 |_{u_0}t))=
$$
$$
=u_0 + 2\log(\cos( |v_0 |_{u_0}at) + \frac{\sgn(a)v_0}{2|v_0
|_{u_0}}\sgn(a)\sin( |v_0 |_{u_0}at)) =
$$
$$
=u_0 + 2\log(\cos( |v_0 |_{u_0}at) + \frac{v_0}{2|v_0
|_{u_0}}\sin( |v_0 |_{u_0}at))= \gamma( at , u_0 , v_0 ),
$$
as claimed, where a property of the functions sinus and cosine is
used.
\end{pf}

\noindent The following result shows an important property shared
by our and the standard exponential maps.

\begin{Lem}\label{lemmadiffenentialexpientity}
Fix a point $u \in \Cb$. The differential of the exponential map
at the origin is the identity operator.
\end{Lem}

\smallskip

\noindent \begin{pf} Let $v \in T_{u} \Cb$ be any point in the
domain of the exponential map. By definition the differential is
$$
(d\exp)_0 v := \frac{d}{dt}_{|t=0} \exp (tv) =\frac{d}{dt}_{|t=0}
\gamma( 1 , u , tv ) =\frac{d}{dt}_{|t=0} \gamma (t , u , v) = v,
$$
where we used the Lemma \ref{lemmahomogeneity}.\end{pf}

\begin{Com}
\emph{In that \\case, as proved by S. Semmes \cite{Sem} and S.
Donaldson \cite{Don}, the geodesic equation is equivalent to a
homogeneous complex Monge-Amp\`ere equation on $M\times
[0,1]\times S^1$. About the Cauchy problem, S. Donaldson in
\cite{Don} constructed examples of short time existence for the
Cauchy problem which cannot be extended and also examples where
the Cauchy problem for geodesics does not admit solution even for
short time. As a consequence, a definition of exponential map is
not at hand as in our case.}
\end{Com}

\section{The Dirichlet problem for the geodesic\\ equation}

In this section we prove that the Dirichlet problem has explicit
unique smooth solution. We think that this was known to Eugenio
Calabi. Moreover, the space $\Cb$ endowed with the Calabi's metric
is a genuine metric space, with a smooth distance function. The
geodesic lines realize the minimum of the length among smooth
paths joining two fixed endpoints. To conclude, we compute that
the diameter of the space is $\frac{\pi}{2}R$.

\begin{Th}\label{teogeodesicallyconvex}
For any two points $u_0 \neq u_1 \in \Cb$ the geodesic line
explicitly given by
$$
e^{\frac{u(t)}{2}} = e^{\frac{u_0}{2}}(\cos(t)+ e^{\frac{u_1 -
u_0}{2}}\frac{\sin(t)}{\sin(t_0)} -
\cos(t_0)\frac{\sin(t)}{\sin(t_0)}),
$$
is the unique (modulo parametrization) solution of the Dirichlet
problem
$$
\left\{
\begin{array}{c}
  D_t u' =0 \\
  u(0) = u_0 \\
  u(t_o) = u_1 ,\\
\end{array}
\right.
$$
where $t_0 \in (0,\frac{\pi}{2})$ is the unique solution of $
\int_M e^{\frac{u_0 + u_1}{2}} \dmuz = \frac{1}{4}\cos(t_0). $
\end{Th}

\smallskip

\noindent \begin{pf} The first claim is that the smooth curve in
the statement is indeed a geodesic. Its velocity, at $t=0$, is
$$
v_0 = \frac{2}{\sin (t_0)} (e^{\frac{u_1 - u_0}{2}}- \cos(t_0));
$$
moreover check that $v_0$ is an element of the tangent space to
$\Cb$ at $u_0$;
$$
\int_M v_0 e^{u_0} \dmuz = \int_M \frac{2}{\sin (t_0)}
(e^{\frac{u_1 - u_0}{2}}- \cos(t_0)) e^{u_0} \dmuz =
$$
$$
=  \frac{2}{\sin (t_0)} \int_M (e^{\frac{u_1 + u_0}{2}}- \cos(t_0)
e^{u_0}) \dmuz = \frac{2}{\sin (t_0)} (\int_M (e^{\frac{u_1 +
u_0}{2}}) \dmuz - \frac{1}{4}\cos(t_0)) =0.
$$
By the uniqueness of solutions for the Cauchy problem the claim
follows. About $t_0$, employ the H\"{o}lder inequality to get
$$
0<\int_M e^{\frac{u_0 + u_1}{2}} \dmuz  \leq (\int_M (e^{\frac{u_0
}{2}})^2 \dmuz)^{\frac{1}{2}}(\int_M (e^{\frac{u_1 }{2}})^2
\dmuz)^{\frac{1}{2}}= \unmezzo \cdot \unmezzo =\frac{1}{4}
$$
and the last inequality is an equality if and only if there exists
a $\lambda \in \Rm$ for which $e^{u_1} = \lambda e^{u_0}$; but the
constraint $\frac{1}{4} = \int_M e^{u_1} \dmuz = \int_M e^{u_0}
\dmuz $ forces $\lambda$ to be $1$, that is $u_0 = u_1$. When
$t\in (0 , \frac{\pi}{2})$, the function $\frac{1}{4}\cos(t)$ is
decreasing and surjective onto $(0, \frac{1}{4})$. Thus there
exists a unique $t_0 \in (0 ,\frac{\pi}{2})$ such that
$$
\frac{1}{4}\cos(t_0) = \int_M e^{\frac{u_0 + u_1}{2}} \dmuz .
$$
Clearly the geodesic is smooth; about the uniqueness, suppose
there are two geodesics connecting $u_0$ and $u_1$. Since we are
interested in uniqueness modulo parametrization, normalize both
geodesics into unitary ones. If $v_0 , w_0 \in T_{u_0} \Cb$ are
the velocity of these two geodesics in $u_0$ and $s,t\in
(-\frac{\pi}{2},\frac{\pi}{2})$ are values for which the point
$u_1$ is reached respectively by the first and by the second
geodesic, i.e.
$$
u_1 = u_0 + 2\log (\cos(s) + \frac{v_0}{2}\sin(s))
$$
and
$$
u_1 = u_0 + 2\log (\cos (t) + \frac{w_0}{2}\sin (t)),
$$
then consider the difference
$$
0 = 2\log(\frac{\cos(s) + \frac{v_0}{2}\sin(s)}{\cos (t) +
\frac{w_0}{2}\sin (t)})
$$
or equivalently
$$
\cos(s) + \frac{v_0}{2}\sin(s)=\cos (t) + \frac{w_0}{2}\sin(t).
$$
Multiply by $e^{u_0}$ to get
$\frac{1}{4}\cos(t)=\frac{1}{4}\cos(s)$ which implies $t=\pm s$
and $\frac{v_0}{2}= \pm \frac{w_0}{2}$; so, if $w_0 = v_0$the two
geodesics coincide, while  if $w_0 = -v_0$ one geodesic is
parametrized backward respect to the other. The theorem is proved.
\end{pf}

\noindent We introduce the following notation. If $\alpha$ is a
smooth curve in $\Cb$ passing through $u_0 , u_1$, we write the
length of $\alpha$ between $u_0$ and $u_1$ as
$$
\Lb_\alpha (u_0 , u_1),
$$
which does not dependent on the parametrization and is symmetric
in $u_0$ and $u_1$.

\begin{Def}
\emph{Let $d$ be the map explicitly defined by
$$
\begin{array}{ccccc}
  d:& \Cb \times \Cb & \rightarrow & [0, \frac{\pi}{2}) &  \\
  &(u_0 , u_1) & \mapsto & d(u_0 , u_1):= & \arccos (4\int_M e^{\frac{u_0 + u_1}{2}}\dmuz)=\Lb_{u} (u_0 , u_1), \\
\end{array}
$$
where $u$ is the geodesic line connecting $u_0$ and $u_1$, and as
remarked in the Theorem \ref{teogeodesicallyconvex} the argument
of $\arccos$ varies between $0$ and $1$; moreover, the
corresponding determination of the function $\arccos$ will range
between $0$ and $\frac{\pi}{2}$. }
\end{Def}

\noindent The next lemma is useful to prove that the function $d$
is a distance and that the geodesic distance is a minimum among
the length of all smooth paths; we follow an argument by Xiu Xiong
Chen \cite{Chen}.

\smallskip

\begin{Lem}\label{triangularinequalityandminimizingproperty}
If $C$ is a smooth curve $C= \phi (s) : [0,1] \rightarrow \Cb$,
then for every $s\in [0,1]$ the length of the geodesic $t \mapsto
u(s,t)$ from the base point $0\in \Cb$ to $\phi (s)$  is not
greater than the length of the geodesic $u(0,t)$ from the base
point $0$ to $\phi (0)$ plus the length of $C$ between $\phi (0)$
and $\phi (s)$; more precisely,  there holds
$$
\Lb_{u(s , t)} (0 , \phi(s)) \leq \Lb_{u(0 , t)} (0 , \phi(0)) +
\Lb_{\phi} (\phi(0) , \phi (s)).
$$
\end{Lem}

\noindent \begin{pf} By definition  the length of $C$ between
$\phi (0)$ and $\phi (s)$ is
$$
\Lb_{\phi} (\phi(0) , \phi (s)) = \int_0^s \sqrt{<\phi_s (\tau) ,
\phi_s (\tau)>_{\phi (\tau)}} d\tau ,
$$
where $\phi_s$ stands for the derivative respect to the parameter
$s$. Suppose, without loss of generality, that $ u(s,t) : [0,1]
\times [0,1] \rightarrow \Cb $ with $u(s,0) =0$, $u(s,1)=\phi(s)$,
for every $s\in [0,1]$. By definition  the length of $u(s,t)$
between $0$ and $\phi (s)$ is
$$
\Lb_{u(s , t)} (0 , \phi(s)) = \int_0^1 \sqrt{<u_{\tau}(s , \tau)
, u_{\tau}(s , \tau)>_{u(s , \tau)}} d\tau .
$$
Notice that by construction there holds
$$
u_s(s,1)=\phi_s(s), \qquad u_s(s,0)=\phi_s(0) \forall s\in [0,1].
$$
Define $ F(s):=\Lb_{\phi} (\phi(0) , \phi (s)) +\Lb_{u(s , t)} (0
, \phi(s))$, so that the claim is
$$
F(s)=\Lb_{\phi} (\phi(0) , \phi (s)) +\Lb_{u(s , t)} (0 , \phi(s))
\geq \Lb_{u(s , t)} (0 , \phi(0)) = F(0).
$$
To get the above claim, it is enough to show $ \frac{\barr
F}{\barr s} \geq 0$. Compute $ \frac{\barr \Lb_{\phi} (\phi(0) ,
\phi (s))}{\barr s}$;
$$
\frac{\barr \Lb_{\phi} (\phi(0) , \phi (s))}{\barr s}
=\sqrt{<\phi_s (s) , \phi_s (s)>_{\phi (s)}} \geq $$ $$\geq - <
\phi_s (s) , u_t (s , 1)>_{\phi(s)} \cdot (\sqrt{<u_t (s , 1) ,
u_t (s , 1)>_{\phi (s)}})^{-1},
$$
where the Cauchy-Schwartz inequality is used. Compute
$$
\frac{\barr \Lb_{u(s , t)} (0 , \phi(s))}{\barr s}= \int_0^1
\frac{1}{2} \frac{1}{\sqrt{<u_t , u_t >_{u(s,\tau)}}} \cdot 2 <D_s
u_t , u_t >_{u(s,\tau)} d\tau =
$$
$$
= \int_0^1 \frac{1}{\sqrt{<u_t , u_t >_{u(s,\tau)}}}<D_t u_s , u_t
>_{u(s,\tau)} d\tau =
$$
$$
=  \frac{1}{\sqrt{<u_t , u_t >_{u(s,\tau)}}}\int_0^1
[\frac{\barr}{\barr t} <u_s , u_t >_u - <u_s , D_t u_t >_u
] d\tau =
$$
$$
= \frac{1}{\sqrt{<u_t , u_t >_{u(s,\tau)}}} [<u_s (s,1) ,
u_t(s,1)>_{\phi(s)} - <0 , u_t (s , 0)>_0
].
$$
Thus
$$
\frac{\barr F}{\barr s} =\frac{\barr \Lb_{\phi} (\phi(0) , \phi
(s))}{\barr s} +\frac{\barr \Lb_{u(s , t)} (0 , \phi(s))}{\barr
s}\geq
$$
$$
\geq (\sqrt{<u_t (s , 1) , u_t (s , 1)>_{\phi (s)}})^{-1} [- <
\phi_s (s) + u_s (s,1) , u_t(s,1)>_{\phi(s)}] =0,
$$
as claimed, where it is used that $\phi_s (s) =u_s (s,1)$ for all
the values of $s$.
\end{pf}

\begin{Th}\label{disadistance}
The map $d$ is a smooth distance function and geodesics realizes
the absolute minimum of the length over all $C^1$ paths.
\end{Th}

\smallskip

\noindent \begin{pf} The fact that $d$ is smooth and symmetric in
$u_0$ and $u_1$ follows from its explicit expression
$$
d(u_0 , u_1) = \arccos (4\int_M e^{\frac{u_0 + u_1}{2}}\dmuz).
$$
Suppose that $d (u_0, u_1) =0 $; then by definition of $d$ the
argument of $\arccos$ must be $1$. On its turn, as remarked in
\ref{teogeodesicallyconvex}, the equality $ \int_M e^{\frac{u_0 +
u_1}{2}}\dmuz = \frac{1}{4}$ holds if and only if $u_0 = u_1$. So
$d(u_0 , u_1) =0$ implies $u_0 = u_1$. To  prove the triangular
inequality, fix any $u_0 , u_1 , u_2 \in \Cb $ and apply the Lemma
\ref{triangularinequalityandminimizingproperty} to the case when
$s=1$ and $C=\phi$ is a geodesic from $\phi (0) = u_1$ and $\phi
(1) = u_2$. The Lemma gives the formula
$$
\Lb_{u(1 , t)} (0 , u_2) \leq \Lb_{u(0 , t)} (0 , u_1) +
\Lb_{\phi} (u_1 ,u_2).
$$
Since we have the freedom to choose the base point $0\in \Cb$
arbitrarily, with the choice $u_0 = 0$ the triangular inequality
follows. About the minimizing property, for any smooth curve $C=
\phi (s) : [0,1] \rightarrow \Cb$, we want to show that the length
of $C$ between $\phi(0)$ and $\phi(1)$ is greater than the
geodesic distance from $\phi (0)$ to $\phi (1)$. This follows from
the lemma proved above, applied to the case when $s=1$ and the
base point $0$ is in $\phi (s)$; in fact,
$$
d(\phi (1) , \phi (0)) \leq d(\phi (1) , \phi (1)) + \Lb_{\phi}
(\phi (0) ,\phi (1)) =\Lb_{\phi} (\phi (0) ,\phi (1)),
$$
as wanted.
\end{pf}

\noindent We conclude this section discussing the diameter of the
space $\Cb$ endowed with the Calabi's metric. So far we proved
that any two points $u_0 , u_1 \in \Cb$ have distance in the range
$[0 , \frac{\pi}{2})$. We claim that the supremum is precisely
$\frac{\pi}{2}$.

\begin{Def}
\emph{The diameter of the space of K\"{a}hler metrics endowed with
the Calabi's metric is
$$
\Diam (\Cb , \Ca):= \sup_{u_0 , u_1 \in \Cb} d(u_0 , u_1).
$$}
\end{Def}

\smallskip

\begin{Th}\label{teodiameter}
There exists a sequence $\{ u_k \}\subset \Cb$ such that, for any
$u \in \Cb$, we have
$$
\lim_{k\rightarrow \infty} d(u , u_k) = \frac{\pi}{2}.
$$
In particular, the space $\Cb$ endowed with the Calabi's metric
has diameter
$$
\Diam (\Cb , \Ca) = \frac{\pi}{2}.
$$
\end{Th}

\smallskip

\begin{pf} Fix any point $p \in M$; for any $k\in \mathbb{N}$ fix two balls $B_{r_k} \subset B_{R_k}$ centered  $p$ with radii
respectively $r_k<R_k$ to be determined. With a standard argument
based on the partitions of unity (see for example \cite{War}), get
a real smooth function $f_k$ on $M$ such that
$$
\left\{
\begin{array}{cc}
  f_k(q)=1 & \mbox{ when $q$ belongs to $B_{r_k}$} \\
  0\leq f_k(q) \leq 1 & \mbox{ when $q$ belongs to $B_{R_k} \backslash B_{r_k}$} \\
  f_k(q)=0 & \mbox{ when $q$ belongs to $M \backslash B_{R_k}$}. \\
\end{array}
\right.
$$
Define
$$
\alpha_k := \int_M f_k \frac{\omega^n}{n!}.
$$
Of course $|B_r| \leq \alpha \leq |B_R|$, where  $|A|$ is the
measure of the open set $A\subset M$. Let $\{ \epsilon_k\}$ be
some positive sequence which tends to zero; the sequence $u_k$ is
defined by
$$
e^{u_k}:= \frac{\frac{f_k}{4\alpha_k} + \epsilon_k }{1 +
\epsilon_k }.
$$
For any $k$, $u_k$ is smooth and satisfies $\int_M e^{u_k}
\frac{\omega^n}{n!} $$ =\frac{1}{4} = Vol(M) $. Fix the radii such
that $ |B_{r_k}| = \frac{1}{4^k}, \; |B_{R_k}| = \frac{1}{3^k}$.
Thus we have
$$
0<\int_M e^{\frac{u_k}{2}} \frac{\omega^n}{n!} \leq
\frac{1}{\sqrt{1+\epsilon_k}} [|M\backslash B_{R_k}|
\sqrt{\epsilon_k} + \frac{|B_{R_k}|}{\sqrt{\alpha_k}}
\sqrt{\frac{1}{4} + \epsilon_k \alpha_k}] \leq
$$
$$
\leq \frac{1}{\sqrt{1+\epsilon_k}} [\frac{1}{4} \sqrt{\epsilon_k}
+ \left(\frac{2}{3}\right)^k \sqrt{\frac{1}{4} + \epsilon_k
\alpha_k}] \rightarrow 0,
$$
when $k\rightarrow \infty$. So if $u$ is any point of $\Cb$,
$$
d(u  , u_k ) = \arccos (4\int_M e^{\frac{u +
u_k}{2}}\frac{\omega^n}{n!}) \geq
$$
$$
\geq \arccos (4 \max_M (u) \int_M
e^{\frac{u_k}{2}}\frac{\omega^n}{n!}) \rightarrow \arccos (0) =
\frac{\pi}{2},
$$
which proves the statement.
\end{pf}

\smallskip

\begin{Rem}
\emph{If the volume of $M$ and thus the sectional curvature $R$ is
not normalized, in general we have
$$
\Diam(\Cb , \Ca) = \frac{\pi}{2}R.
$$}
\end{Rem}

\smallskip

\begin{Def}\emph{
We call the boundary of $\Cb$ and we denote it by $\barr\Cb$, the
space
$$
\barr\Cb := \{ u\in C^{\infty} (M , \Rm \cup \{ - \infty\} ) \, |
\, \int_M e^u \dmuz = \int_M \dmuz , \; \exists x\in M : \, u(x) =
-\infty\}.
$$
$\square$ }\end{Def}

\smallskip

\begin{Th}\label{teodistancefromtheboundaryiszero}
For any fixed $u_0 \in \Cb$ there exists a sequence $\{ w_k \}
\subset \barr \Cb$ such that
$$
d(u_0 , w_k) \rightarrow 0.
$$
In particular the distance of $u_0 \in \Cb$ from the boundary
$\barr\Cb$ is zero.
\end{Th}

\smallskip

\noindent \begin{pf} A unitary geodesic starting from $u_0$ with
initial velocity $v_k\in T_{u_0} \Cb$ has explicit expression $
e^{\frac{u(t)}{2}}= e^{\frac{u_0}{2}}(\cos(t) +
\frac{v_k}{2}\sin(t)). $ The parameter $t$ corresponds to the
distance between $u_0$ and $u(t)$; moreover, when $ 0\leq t =
\arccot (-\frac{\min(v_k)}{2}) $, then $u(t)$ belongs to the
boundary $\partial \Cb$. With a technique very similar to that one
employed in the proof of Theorem \ref{teodiameter}, is possible to
construct a sequence of unitary vectors $\{v_k \}_k \subset
T_{u_0} \Cb$ such that $ \lim_{k\rightarrow \infty} \min v_k =
-\infty . $ We claim that with such a sequence we are able to
define a sequence of points of the boundary of $\Cb$ such that the
distance rom $u_0$ goes to zero. In fact, define, for all $k\in
\mathbb{N}$,
$$
w_k := u (t_k) \in \barr\Cb
$$
and compute
$
d(u_0 , w_k ) = d(u_0 , u(t_k)) = t_k = \arccot (-\frac{\min
v_k}{2}),
$
so that
$$
\lim_{k\rightarrow \infty} d(u_0 , w_k) = 0,
$$
as requested by the statement.
\end{pf}

\begin{Com}
\emph{In that \\
case the construction of the Chen's geodesics is completely
different from ours. In fact Xiu Xiong Chen \cite{Chen} had to
deal with a degenerate fully non linear second order partial
differential equation on the manifold $M \times S^1 \times [0,1]
$. He used the continuity method; thus he wrote a family of
Monge-Amp\`ere equations, depending on a parameter $\lambda$, such
that at $\lambda =0$ the equation was the geodesic one whose
solvability was required to be shown. The method consisted of
proving that the set of $\lambda \in [0,1]$ whose associated
problem admits a solution is not empty, open and closed. The first
two points were not difficult to prove, since non degenerate
Monge-Amp\`ere on manifold with boundary had already been
understood. The crucial point was about the closedness, i.e. the
proof of a priori estimates for the solutions. The $C^0$ bound
from above was given by the use of the maximum principle, while
the bound from below was got by a trick that goes back to E.
Calabi. The $C^0$ estimate were enough to get interior estimates
from above of the Laplacian; for this Chen employed a lemma by
Yau. The estimates of the Laplacian from below were immediate.
Thus Chen passed to consider the remaining case, i.e. the
estimates form above of the Laplacian at the boundary of the
manifold. About this Chen first proved that the mixed derivatives
with a tangential component dominates the other ones; then, with a
surprising barrier function, Chen showed that those derivatives
are controlled by the gradient norm. The last step was to control
the gradient of the norm from above. This step would give also the
Laplacian estimates, for what we just said. The estimate of the
gradient of a solution was got via a blowing-up analysis argument.
Chen's theorem partially fulfills one conjecture in Donaldson's
program \cite{Don}, since Chen's geodesics are $C^{1,1}$ and it is
not known yet if they are smooth. Moreover, it is still not known
if Chen's geodesics lay completely inside the space, as in our
case. Again by a Xiu Xiong Chen's theorem in \cite{Chen}, the
space of K\"{a}hler metric endowed with that metric is a genuine
metric space, and geodesics minimize the length among smooth
paths. This completely fulfills another conjecture in Donaldson's
program \cite{Don}. We got the same result in the Proposition
\ref{disadistance}, with the difference that in our case the
distance function is smooth and in that case is $C^1$.}
\end{Com}

\section{The Jacobi equation; $\Cb$ as a portion of a sphere}
\noindent In this section we use the nice expression of geodesic
lines and sectional curvature to define Jacobi fields and
conjugate points in the same way of Riemannian geometry. We give a
characterization and a link between conjugate points and the
exponential maps; we prove that there are no conjugate points.
Finally we show an isometric immersion of $\Cb$ in a portion of
the sphere in $C^{\infty}(M, \Rm)$.

\begin{Def}
Let $u = u(t) : (-\epsilon , \epsilon ) \rightarrow \Cb$ be  a
geodesic, and let $J= J(t)$ be a smooth section along $u$. We call
$J$ a Jacobi field if it satisfies the Jacobi equation
$$
D_t^2 J - R(u ' , J ) u' =0.
$$
\end{Def}

\begin{Rem}
\emph{Fix a point $u_0 \in \Cb$ in the space of K\"{a}hler
metrics. Suppose to have a curve $v= v(s)$ in $\Vb$, and consider
the two parameter family $(s,t) \mapsto \exp (tv(s))$, where we
wrote $\exp$ instead $\exp_{u_0}$. Notice that, when $s$ is fixed,
the curve $ \alpha(t):=\exp(tv(s))$ is a geodesic for any fixed
$s$. In fact, by the Lemma \ref{lemmahomogeneity} this curve can
be expressed by
$$
t\mapsto \gamma (1 , u_0 , tv(s)) = \gamma (t , u_0 , v(s))
$$
which is precisely the form of a geodesic, by the definition of
$\gamma$. We a smooth section $J$ along the curve $\alpha$ as
$$
J (t):=  \mapsto \frac{\barr}{\barr s}\exp (tv(s)).
$$
This is an example of a Jacobi field as we compute
$$
D_t^2 \frac{\barr}{\barr s}\exp (tv(s)) = D_t D_s
\frac{\barr}{\barr t} \exp (tv(s)) =
$$
$$
= (D_t D_s - D_s D_t )\frac{\barr}{\barr t} \exp (tv(s))= R(\alpha
' , J ) \alpha ',
$$
where we used in the last equality that $\alpha$ satisfies the
geodesic equation.}
\end{Rem}

\begin{Lem}\label{lemmaontwoformulasapriorijacobiù}
For a Jacobi field $J$ along a geodesic $u$ these two formulae
hold
$$
<u' , D_t J >(t)=<u' , D_t J >(0), \quad <u' ,  J >(t) = <u' , D_t
J >(0) \cdot t + <u' , J >(0).
$$
\end{Lem}

\noindent \begin{pf} First notice that
$$
\frac{\barr}{\barr t} <u' , D_t J > = <u' , D_t^2 J > = <u' , R(u
' , J ) u'>=
$$
$$
=<u' , <u' , u' > J - <u' ,  J> u'>= 0.
$$
Thus $ <u' , D_t J > (t) =  <u' , D_t J > (0)$. About the second
formula,
$$
\frac{\barr}{\barr t} <u' , J > = <u' , D_t J>(t) = <u' , D_t
J>(0),
$$
which gives, after an integration, the second claim.\end{pf}

\begin{Lem}\label{lemmasecondorfereodejacobi}
The Jacobi equation is equivalent to the second order differential
equation
$$
J'' + u' J' -2|v_0|_{u_0}^2 J +4<v_0 , D_t J (0)>_{u_0} + 2 (<v_0
, D_t J (0)>_{u_0} t + <v_0 , J(0)>_{u_0})=0,
$$
where $u(t)$ is a geodesic with initial data $(u_0 , v_0)$ and $J$
is a Jacobi field along $u$.
\end{Lem}

\smallskip

\noindent \begin{pf} In the next computation is assumed $R=1$;
$$
D_t^2 J = D_t (J ' + \frac{1}{2}u' J + 2 <u' , J >) = J'' +
\frac{1}{2} u' J' + 2 <u' , J' >+\frac{1}{2} u'' J +
$$
$$
+\frac{1}{2} u' J'  +\frac{1}{4} (u')^2 J +
 <u' , u' J> +2<u' , J > + <u' , J >u' =
$$
$$
=J'' + u' J' + 2<u' , J' + \unmezzo u' J  > + \unmezzo J (u'' +
\unmezzo (u')^2) + 2 <u' , D_t J> + <u' , J>u' =
$$
$$
=J'' + u' J' + 4<u' , D_t J>+ \unmezzo J (u'' + \unmezzo (u')^2)+
<u' , J>u'
$$
$$
=J'' + u' J' + 4<u' , D_t J>-  <u' , u'>J +<u' , J>u'.
$$
Thus the Jacobi equation is equivalent to
$$
D_t^2 J - R(u' , J) u' =
$$
$$
=J'' + u' J' + 4<u' , D_t J>- <u' , u'>J + <u' , J>u' -<u' , u'> J
+  <u' , J> u'=
$$
$$
=J'' + u' J' + 4<u' , D_t J>- 2<u' , u'>J + 2<u' , J>u' =
$$
$$
=J'' + u' J' +4<v_0 , D_t J(0)>- 2 |v_0 |_{u_0}^2 J + 2u' (<v_0 ,
D_t J(0)>t + <v_0 , J(0)>),
$$
where in the last equality the Lemma
\ref{lemmaontwoformulasapriorijacobiù} is applied. The claim is
proved.\end{pf}

\smallskip

\begin{Lem}\label{proponuniquenessofjacobi}
Let $u$ be a geodesic in $\Cb$. If $J$ and $\tilde{J}$ are Jacobi
fields such that $J(0) = \tilde{J}(0)$ and $D_t J (0) = D_t
\tilde{J} (0)$ then $J = \tilde{J}$.
\end{Lem}

\smallskip

\noindent \begin{pf} By the  Lemma
\ref{lemmasecondorfereodejacobi}, the Jacobi equation is a second
order  differential equation for which the theorem on uniqueness
of solution applies.
\end{pf}

\begin{Prop}\label{propcharacterizationjacobifields}
Let $u = u(t)$ be a geodesic in $\Cb$ with initial velocity $v_0
\in T_{u_0} \Cb$ and let $J$ be a Jacobi field along $u$ with
$J(0)=0$ and $D_t J (0) = w$. Let $v(s)$ be a curve in $T_{u_0}
\Cb$ with $v(0) = v_0$ and $v'(0) = w$. Consider the Jacobi field
$\tilde{J} (t) = \frac{\barr}{\barr s}_{s=0} \exp (tv(s)). $ Then
$J = \tilde{J}$.
\end{Prop}

\smallskip

\noindent \begin{pf} First notice that $w$ is an element of
$T_{u_0} \Cb$ since it is $w = D_t J(0) \in T_{u_0} \Cb$. Also,
$v'(0)$ can be identified with an element of $ T_{u_0} \Cb$; in
fact , apply the operator $\frac{\barr}{\barr s}_{|s=0}$ to the
identity$ \int_M v(s) e^{u_0} \dmuz =0$ and get
$$
\int_M \frac{\barr v}{\barr s}_{|s=0} e^{u_0} \dmuz =0.
$$
It is clear that $\tilde{J} (0) =0$; moreover notice that
$$
D_t \frac{\barr}{\barr s}_{s=0} \exp (tv(s)) = D_t
((d\exp)_{tv_0}(tw) ) =
$$
$$
= D_t (t(d\exp)_{tv_0}(w) ) =(d\exp)_{tv_0}(w) + t D_t
((d\exp)_{tv_0}(w) ).
$$
In particular for $t=0$ the above formula reads
$$
D_t \tilde{J}(0) =(d\exp)_{0}(w) = w,
$$
where the Lemma \ref{lemmadiffenentialexpientity} is used.
Conclude, by the Lemma \ref{proponuniquenessofjacobi}, that $J =
\tilde{J}$.\end{pf}

\noindent We want to show a link between the exponential map and
the Jacobi fields that occurs as well in the finite dimensional
theory. We start with the following notion

\begin{Def}
\emph{ Let $u_0$ and $u_1$ be two points of $(\Cb , \Ca)$, and
$u=u(t)$ the geodesic that joins them so that $u(0) = u_0$ and
$u(t_0) = u_1$. The two points $u_0 , u_1$ are called conjugate
when there is a Jacobi field $J$ along $u$, $J$ non identically
zero, such that $J (0)=0$ and $J(t_0) = 0$.}
\end{Def}

\begin{Prop}
Let $u=u(t)$ be a geodesic with $u(0)=u_0$ and $u'(0)\neq 0$. The
point $u_1 = u(t_0)$ is conjugate to $u_0$ if and only if $v_0 :=
t_0 u'(0)\in T_{u_0} \Cb$ is a critical point of $\exp$.
\end{Prop}

\noindent \begin{pf} With the use of the Proposition
\ref{propcharacterizationjacobifields}, write explicitly $J$ as
$$
J(t) = (d\exp)_{t u'(0)}(tD_t J(0)).
$$
The field $J$ is not identically zero if and only if $D_t J(0)
\neq 0$ as follows from the Lemma \ref{proponuniquenessofjacobi}.
Notice that $J(0) = 0$; thus $u_0$ and $u_1 $ are conjugate if and
only if $ J(t_0) =(d\exp)_{t_0 u'(0)}(t_0 D_t J(0)) =0 $, which
means that $ t_0 u'(0) $ is a critical point of $\exp$. \end{pf}

\noindent We want to prove a statement about the existence of
conjugate points in our case. First we need the following result.

\begin{Lem}
If $u_0 , u_1$ are conjugate points and $J (0) = u_0$, $J(a)
=u_1$, $J$ is not identically zero, necessarily we have
$$
<J , u' >(t) =0,\quad <D_t J , u' >(t) = 0.
$$
In particular the Jacobi equation simplifies into
$$
J'' + u' J' - 2 |v_0 |_{u_0}^2 J =0.
$$
\end{Lem}

\smallskip

\noindent \begin{pf} Consider the second  claim; from the Lemma
\ref{lemmaontwoformulasapriorijacobiù}
$$
<J , u' >(t) = <D_t J (0) , u'(0)> t + <J (0) , u' (0)>.
$$
Use that $J(0)=0$ and evaluate the formula at $t= a$ to get $<J ,
u' >(a) = <D_t J (0) , u'(0)> a$. Since $J(a)=0$ the addendum $
<D_t J (0) , u'(0)>$ vanishes and the formula simplifies into
$$
<J , u' >(t) =  <J (0) , u' (0)>,
$$
which is zero identically since $J(0)=0$. About the first claim,
still from  the Lemma \ref{lemmaontwoformulasapriorijacobiù}, $
<u' , D_t J
>(t)=<u' , D_t J >(0)$, which is already proved to be zero. Plug
in the two formulae in the Lemma \ref{lemmasecondorfereodejacobi}
to get the simplified equation.\end{pf}

\smallskip

\begin{Th}\label{teonoconjpoints}
In the space $(\Cb , \Ca)$ there are no conjugate points.
\end{Th}

\smallskip

\noindent \begin{pf} Suppose by contradiction that $u_0 , u_1 \in
\Cb $ are conjugate points. So there exists a not identically zero
Jacobi field $J$ along the geodesic $u$ which connects $u_0 =
u(0)$ and $u_1 = u(a)$. $J$ necessarily satisfies the equation
$$
J'' + u' J' -2|v_0 |_{u_0}^2 J =0
$$
Suppose, without loss of generality, that the geodesic which
connects $u_0$ and $u_1$ is unitary. Thus the Jacobi equation
simplifies into
$$
J'' + u' J' -2 J =0.
$$
Consider
$$
(e^{\frac{u}{2}}J)'' = e^{\frac{u}{2}}J'' + e^{\frac{u}{2}} J' u'
+ e^{\frac{u}{2}}J (\frac{u''}{2} + \frac{(u')^2}{4}) =
e^{\frac{u}{2}} J'' + e^{\frac{u}{2}} u' J' - e^{\frac{u}{2}} J.
$$
Thus the equation necessarily satisfied by $J$ is equivalent to
the following one
$$
(e^{\frac{u}{2}}J)'' - e^{\frac{u}{2}}J =0
$$
and the generic solution is
$$
e^{\frac{u}{2}}J = A \cos (t) + B \sin (t), \quad A , B \in \Rm .
$$
The conditions $J(0)=J(a) =0$ and $J\neq 0$ entail
$$
A=0, \quad \sin(a) =0,
$$
that is $a = k\pi$, where $k \in \mathbb{Z}$, $k\neq 0$. But the
geodesic $u$, exists at most in the interval $(-\frac{\pi}{2},
\frac{\pi}{2})$, so a number $a= k\pi$, with $k \in \mathbb{Z}$,
$k\neq 0$ cannot be contained in the interval of definition of
$u$. This contradicts the fact that the Jacobi field $J$ exists
for $t=a$.\end{pf}

\smallskip

\begin{Rem}
\emph{The Theorem \ref{teonoconjpoints} is based on the fact that
in the geometry arisen from the Calabi's metric, geodesic are not
long enough to develop conjugate points. In Riemannian geometry a
conjugate point is also linked to the point where a geodesic
ceases to be a minimum of the length among smooth paths. Thus the
above theorem is somehow consistent with the Proposition
\ref{disadistance}.}
\end{Rem}

\smallskip

\noindent To conclude this section, we show an isometric immersion
of $\Cb$ with the Calabi's metric into a portion of the sphere of
$C^{\infty}(M , \Rm)$ with a flat metric. Consider the space
$C^{\infty}(M, {\mathbb R})$ endowed with the metric
$$
\prec \psi , \chi \succ_\phi := \int_M \psi \chi
\frac{\omega^n}{n!}, \quad \phi \in C^{\infty}(M  ,{\mathbb R}),
\, \psi , \chi \in T_\phi C^{\infty}(M , {\mathbb R}),
$$
which is flat since there is not dependence on the point $\phi$
where the scalar product is computed. We call $\prec \cdot , \cdot
\succ$ the Euclidean metric on $C^{\infty}(M , {\mathbb R})$.

\begin{Th}\label{teoimmersionsphere}
The map $ A : {\mathcal C} \rightarrow C^{\infty}(M , {\mathbb R})
$ defined as $A(u) = 2e^{\frac{u}{2}}$ is injective and has image
$ A(\Cb)= \{ f\in C^{\infty} (M , {\mathbb R}) \, | \, f>0 ,\,
\int_M f^2 \frac{\omega^n}{n!} = 1 \}. $ Moreover, the pullback of
the Euclidean metric on $C^{\infty}(M , {\mathbb R})$ via the map
$A$ is the Calabi's metric.
\end{Th}

\smallskip

\begin{pf}
\emph{The differential of $A$ is
$$
(dA)_u (v) = \frac{d}{dt}_{|t=0} A(\exp_u (tv)) =
2e^{\frac{u}{2}}\frac{1}{2} (d\exp_u)_0 (v) = e^{\frac{u}{2}} v .
$$
Thus the pullback is
$$
(A^{*}Eucl)[u](v,w) = \prec (dA)_u (v), (dA)_u(w) \succ_{A(u)} =
\int_M e^{\frac{u}{2}}v e^{\frac{u}{2}} w \frac{\omega^n}{n!} =
<v, w>_u ,
$$
as claimed.}
\end{pf}

\begin{Rem}
\emph{ The map $A$ is an isometric inclusion in analogy with the
case of the isometric immersion of the part of the sphere $S^2$
lying in the first quadrant;
$$
\{ x\in {\mathbb R}^3 \, | \, x_1 >0 , x_2 >0 , x_3
>0 , \sum_{k=1}^3 x_k^2 =1\} \rightarrow ({\mathbb R}^3 , Eucl).
$$
Moreover, the diameter of this $2$ dimensional pece of the sphere
is $\frac{\pi}{2}$, consistently with the formula of the diameter.
Everything scales correctly when the radius of the sphere is $R$
instead of $1$. Namely we have, when the radius of ${\mathcal C}$
is $R=\frac{1}{2\sqrt{Vol}}$ with a generic value of the volume of
$M$, that the map $A$ sends ${\mathcal C}$ in
$$
\{ f\in C^{\infty}(M , {\mathbb R}) \, | \, f>0, \, \int_M f^2 =
R^2\};
$$
in this case the sectional curvature  of ${\mathcal C}$ is
$\frac{1}{R^2}$ and the diameter   is $\frac{\pi}{2}R$.}
\end{Rem}

\section{The Calabi's gradient metric}\label{sectiongradientnorm}

\smallskip

\noindent Again from  an idea by Calabi \cite{CC2} comes the
metric
$$
\ll \psi , \chi  \gg_\phi := \int_M (\nabla \psi , \nabla \chi
)_\phi \dmu
$$
that we call the Calabi's gradient metric. It is a well defined
metric on the space  $\tilde{\Hb}$. We are going to show the
existence of the Levi Civita covariant derivative for the Calabi's
gradient metric. Then we are going to discuss the curvature
induced by the Calabi's gradient metric and the geodesic equation
in the case when $M$ is a closed Riemann surface. In the proof of
the existence of the Levi Civita covariant derivative we will use
the following basic result of Riemannian geometry due to W. Hodge,
which we state without a proof.

\begin{Lem}
Let $(X , g)$ be a closed Riemannian manifold of dimension $n$,
and let $\Delta_g$ the Laplacian metric operator of the given
Riemannian metric $g$. If a smooth real function $ f\in C^{\infty}
(X , \Rm) $ satisfies
$$
\int_M f  d \mu_g =0
$$
then there exists a smooth real function $h\in C^{\infty}(X ,
\Rm)$ such that
$$
f = \Delta_g h .\; \Box
$$
\end{Lem}

\smallskip

\begin{Prop}
The Levi Civita covariant derivative for the Calabi's gradient
metric is given by
$$
2\Delta_\phi D_t \psi = \ddt (\Delta_\phi \psi ) + \Delta_\phi
\psi ' + (\Delta_\phi \psi )(\Delta_\phi \phi '),
$$
where $\phi$ is a smooth curve in $\tilde{\Hb}$ and $\psi$ is  a
smooth section along $\phi$.
\end{Prop}

\smallskip

\noindent \begin{pf} We notice that the Calabi's gradient metric
can be written, integrating by parts, as
$$
\ll \psi ,\chi \gg_\phi = -\int_M \psi \Delta_\phi \chi \dmu =
-\int_M \chi \Delta_\phi \psi \dmu .
$$
Thus, if $\phi = \phi(t)$ is a smooth path in the space of
K\"{a}hler metrics and $\psi = \psi (t)$ is a smooth real  section
on $\phi$ we compute, to get the compatibility with the metric,
$$
\ddt (- \int_M \psi \Delta_\phi \psi \dmu )= -\int_M (\psi '
\Delta_\phi \psi + \psi \ddt(\Delta_\phi \psi ) + \psi
(\Delta_\phi \psi )(\Delta_\phi \phi '))\dmu =
$$
$$
= - \int_M \psi ( \Delta_\phi \psi ' + \ddt(\Delta_\phi \psi ) +
(\Delta_\phi \phi ' )(\Delta_\phi \psi))\dmu .
$$
The computation here above suggests to define implicitly the Levi
Civita covariant derivative as
$$
2\Delta_\phi D_t \psi = \ddt (\Delta_\phi \psi ) + \Delta_\phi
\psi ' + (\Delta_\phi \psi )(\Delta_\phi \phi ').
$$
The expression on the righthand side has integral zero; in fact
$$
\int_M  \ddt (\Delta_\phi \psi ) + \Delta_\phi \psi ' +
(\Delta_\phi \psi )(\Delta_\phi \phi ')\frac{\omega_\phi^n}{n!} =
\frac{d}{dt}\int_M \Delta_\phi \psi \frac{\omega_\phi^n}{n!} +
\int_M \Delta_\phi \psi ' \frac{\omega_\phi^n}{n!} = 0,
$$
where is just used integration by parts. Apply the Hodge result to
get that there exists a function $A_{\phi} : T_{\phi} \tilde{\Hb}
\rightarrow C^{\infty}(M , \Rm)$ such that
$$
\Delta_{\phi} A =\ddt (\Delta_\phi \psi ) + \Delta_\phi \psi ' +
(\Delta_\phi \psi )(\Delta_\phi \phi ').
$$
Now define $ D_t : T_{\phi} \tilde{\Hb} \rightarrow T_{\phi}
\tilde{\Hb} $ this way
$$
D_t \psi :=A_{\phi}(\psi) - L(0, A_{\phi}(\psi)).
$$
Finally, it is straightforward to check the other properties of a
Levi Civita covariant derivative.
\end{pf}

\smallskip

\noindent When $M$ is a closed Riemann surface, the Levi Civita
covariant derivative has the explicit formulation
$$
D_t \psi = \psi ' - \frac{1}{\Vol}\ll \psi , \phi ' \gg_\phi .
$$
This makes possible an easy computation of the curvature tensor
and of the geodesic curves

\smallskip

\begin{Prop}
When $M$ is a closed Riemann surface, the space of K\"{a}hler
metrics endowed with the Calabi's gradient metric has zero
curvature tensor and hence  zero sectional curvature.
\end{Prop}

\smallskip

\noindent \begin{pf} The second claim is a direct consequence of
the first claim. To prove the first claim, we just remark that, if
$\phi (s,t)$ is a smooth two parameter family in the space of
K\"{a}hler metrics and $\psi$ is a section along it, then
$$
D_t D_s \psi =D_t( \psi_s - \frac{1}{\Vol}\ll \psi , \phi_s
\gg_\phi) = \psi_{ts} - \frac{1}{\Vol}\ll \psi_s , \phi_t \gg_\phi
-\ddt \frac{1}{\Vol}\ll \psi , \phi_s \gg_\phi =
$$
$$
=\psi_{ts} - \frac{1}{\Vol}(\ll \psi_s , \phi_t \gg_\phi + \ll
\psi_t , \phi_s \gg_\phi)- \frac{1}{\Vol}\ll \psi , \phi_{ts}
\gg_\phi
$$
which is symmetric on $s$ and $t$, that is the curvature tensor is
zero as claimed.
\end{pf}

\begin{Th}
When $M$ is a closed Riemann surface, the solution of the Cauchy
problem
$$
\left\{
\begin{array}{l}
  D_t \phi ' =0 \\
  \phi (0)=\phi_0 \\
  \phi ' (0) = \phi_0 '\\
\end{array}
\right.
$$
is
$$
\phi (t) = \frac{1}{2Vol}\ll \phi_0 ' , \phi_0 ' \gg_{\phi_0} t^2
+ \phi_0 ' t + \phi_0 .
$$
\end{Th}

\smallskip

\noindent \begin{pf} Clearly the curve
$$
\phi (t) = \frac{1}{2Vol}\ll \phi_0 ' , \phi_0 ' \gg_{\phi_0} t^2
+ \phi_0 ' t + \phi_0
$$
solves the differential problem; thus it remains only to check
that the curve lies inside $\tilde{\Hb}$, that is $\int_M \phi'(t)
\omega_{\phi(t)} =0$. In fact
$$
\int_M \phi'(t) \omega_{\phi(t)} = \int_M (\frac{1}{Vol}\ll \phi_0
' , \phi_0 ' \gg_{\phi_0} t + \phi_0 ')(1 + t\Delta_{\phi_0}
\phi_0 ' )\omega_{\phi_0} =
$$
$$
=\ll \phi_0 ' , \phi_0 ' \gg_{\phi_0} t + \int_M \phi_0
'\omega_{\phi_0} +  \frac{1}{Vol}\ll \phi_0 ' , \phi_0 '
\gg_{\phi_0} t^2  \int_M \Delta_{\phi_0} \phi_0 ' \omega_{\phi_0}+
\int_M t\phi_0 ' \Delta_{\phi_0} \phi_0 ' \omega_{\phi_0}=
$$
$$
= t(\ll \phi_0 ' , \phi_0 ' \gg_{\phi_0} - \int_M t\phi_0 '
\Delta_{\phi_0} \phi_0 ' \omega_{\phi_0}) =0,
$$
where is used that $ \int_M \phi_0 '\omega_{\phi_0}=0 $ and that
clearly $ \int_M \Delta_{\phi_0} \phi_0 ' \omega_{\phi_0} =0$.
Notice that the addendum of second degree in $t$ does not depend
on space variables, thus at the level of K\"{a}hler metrics the
expression of the geodesic is affine  in $t$ as expected when the
curvature is zero.
\end{pf}

\smallskip

 \vskip3cm

 \vskip2cm
\noindent
Simone Calamai\\
Dip. Matematica ``U. Dini''\\
Universit\`a di Firenze\\
Viale Morgani 67/a\\
I-50134 Firenze -- ITALY


\begin{thebibliography}{20}


\bibitem{AT}
C.\ Arezzo and G.\ Tian, \textsl{Infinite geodesic rays in the
space of K\"{a}hler \\ potentials}, Ann. Sc. Norm. Super. Pisa Cl.
Sci. (5) 2 (2003), no. 4, 617--630.

\bibitem{Bes}
A.\ Besse, \textsl{Einstein manifolds}, Classics in Mathematics.
Berlin: Springer, 1987.

\bibitem{B}
J.P.\ Bourguignon, \textsl{Ricci curvature and measures}, Japan.
J. Math. 4, 27-45 (2009).

\bibitem{Cal0}
E.\ Calabi, \textsl{The space of K\"{a}hler metrics}, , Proc.
Internat. Congr. Math. (Amsterdam, 1954), Vol 2, Noordhoff,
Groningen, and North-Holland, Amsterdam, 1956, pp. 206-207.

\bibitem{Cal}
E.\ Calabi, \textsl{Extremal K\"{a}hler metrics}, In:\ "Seminar on
Differential\\ Geometry", Ann. Math. Stud. \ {\bf 102} (1982),
259--290.

\bibitem{CalII}
E.\ Calabi, \textsl{Extremal K\"{a}hler metrics II}, In:\
"Differential Geometry and \\ Complex analysis", Lecture notes in
Math., Springer, 1985, 96--114.

\bibitem{CC}
E.\ Calabi and X.X.\ Chen, \textsl{The space of K\"{a}hler metrics
II}, J. Differential Geom.  \ {\bf 56} (2000), 189--234.

\bibitem{CC2}
Private communication between X.X.\ Chen and E.\ Calabi.

\bibitem{Chen}
X.X.\ Chen, \textsl{The space of K\"{a}hler metrics}, J.
Differential Geom.  \ {\bf 61} (2002), 173--193.

\bibitem{ChenIII}
X.X.\ Chen, \textsl{Space of Kähler metrics, III. On the lower
bound of the Calabi energy and geodesic distance.},  Invent. Math.
\ {\bf 175} (2009), no. 3, 453--503.

\bibitem{ChenIV}
X.X. Chen, \textsl{Space of Kähler metrics, IV. On the lower bound
of the K-energy.} Preprint.

\bibitem{CH}
X.X.\ Chen, W.\ He, \textsl{The space of volume forms} Preprint.

\bibitem{ChenSun}
X.X. \ Chen; S. \ Sun, \textsl{Space of K\"{a}hler metrics (V).
K\"{a}hler quantization.} Preprint.


\bibitem{CT}
X.X.\ Chen; G.\ Tian , \textsl{Geometry of Kahler metrics and
foliations by holomorphic discs }, Publ. Math. Inst. Hautes Etudes
Sci. No. 107 (2008), 1--107.

\bibitem{Don}
S.K.\ Donaldson, \textsl{Symmetric spaces, K\"{a}hler geometry and
Hamiltonian \\dynamics}, In:\ "Northern California Symplectic
Geometry Seminar",\ 13--33,\ Amer.\ Math.\ Soc.\ Transl.\ (2) \
{\bf 196} Amer.Math. Soc., Providence, RI , 1999.

\bibitem{Don2}
S.K.\ Donaldson, \textsl{Nahm's equation and free boundary
problems}, Preprint.

\bibitem{Mab}
T.\ Mabuchi, \textsl{Some symplectic geometry on compact
K\"{a}hler manifolds}, \\Osaka J. Math. \ {\bf 24} (1999),
227--252.

\bibitem{Mab2}
T.\ Mabuchi, \textsl{ $K$-energy maps integrating Futaki
invariants}, \\Tohoku Math. J. (2) \  {\bf 38}  (1986),  no. 4,
575--593.

\bibitem{Mil}
J.\ Milnor, \textsl{ Morse Theory }, Princeton Uniersity Press,
1963.

\bibitem{Sem}
S.\ Semmes, \textsl{Complex Monge-Amp\`ere and symplectic
manifolds}, Amer. J. Math., \ {\bf 114} (1999), 495--550.

\bibitem{War}
F. W.\ Warner, \textsl{Foundations of Differentiable Manifolds and
Lie Groups}, Springer-Verlag, 1983.

\bibitem{Yau}
S. T.\ Yau, \textsl{On the Ricci Curvature of a Compact K\"{a}hler
Manifold and the Complex Monge-Amp\`ere Equation}, Comm. P. And
Appl. Math., 31(1978), 339-411.

\end{thebibliography}
\end{document}